\documentclass[11pt]{article}

\usepackage{amssymb,amsmath,amsthm,amsfonts}

\usepackage{epsfig}

\usepackage{rotating}

\usepackage{subfigure}

\newcommand{\pd}[2]{\frac{\partial#1}{\partial#2}}

\begin{document}

\title{Mori-Zwanzig reduced models for uncertainty quantification I: Parametric uncertainty}
\author{Panos Stinis \\ 
Department of Mathematics \\
University of Minnesota \\
    Minneapolis, MN 55455} 
\date {}

\maketitle

\begin{abstract}
In many time-dependent problems of practical interest the parameters entering the equations describing the evolution of the various quantities exhibit uncertainty. One way to address the problem of how this uncertainty impacts the solution is to expand the solution using polynomial chaos expansions and obtain a system of differential equations for the evolution of the expansion coefficients. We present an application of the Mori-Zwanzig formalism to the problem of constructing reduced models of such systems of differential equations. In particular, we construct reduced models for a subset of the polynomial chaos expansion coefficients that are needed for a full description of the uncertainty caused by the uncertain parameters. The viscous Burgers equation with uncertain viscosity parameter is used to illustrate the construction. For this example we provide a way to estimate the necessary parameters that appear in the reduced model {\it without} having to solve the full system.
\end{abstract}

\section{Introduction}

The problem of quantifying the uncertainty of the solution of systems of partial or ordinary differential equations has become in recent years a rather active area of research. The realization that more often than not, for problems of practical interest, one is not able to determine the parameters, initial conditions, boundary conditions etc. to within high enough accuracy, has led to a flourishing literature of methods for quantifying the impact that this uncertainty imposes on the solution of the problems under investigation (see e.g. \cite{ghanem,leonenko,ma,nouy,venturi,wan}). However, despite the increase in computational power and the development of various techniques for uncertainty quantification there is still a wealth of problems where reliable uncertainty quantification is beyond reach. The main reason behind the inadequacy is the often high dimensionality (in probability space) of the uncertainty sources. When this uncertainty is coupled with the fact that for practical problems, even solving the corresponding equations for one value of the uncertain parameter (initial condition, boundary condition, $\ldots$) can be very expensive, it results in the uncertainty quantification problem being a rather formidable task. One way to address this problem is to look for reduced models for a subset of the variables needed for a complete description of the uncertainty. 

We begin by noting that not all sources of uncertainty are created equal. For example, as we are taught from the theory of ordinary and partial differential equations, the effect of uncertainty in the initial conditions is different from the effect of a parametric uncertainty (see e.g. \cite{barreira,evans}). In addition, the effect of all types of uncertainty is intimately connected with the inherent instabilities that may be present in the underlying system which we subject to the uncertainty. These considerations remain equally, if not more, important when we attempt to construct reduced models for uncertainty quantification. 

In the current work, we are concerned with the construction of reduced models for systems of differential equations that arise from polynomial chaos expansions of solutions of a PDE or ODE system. In particular, we focus on the case that the given PDE or ODE system contains some uncertain parameter and we want to construct a reduced model for the evolution of a subset of the polynomial chaos expansions that are needed for a complete description of the uncertainty caused by the uncertain parameters. There are different methods to construct reduced models for PDE or ODE systems (see e.g. \cite{givon,CS05} and references therein). We choose to use the Mori-Zwanzig (MZ) formalism in order to construct the reduced model \cite{CHK00,CHK3}. 

The main issue with all model reduction approaches is the computation of the memory caused by the process of eliminating variables from the given system (referred to as the full system from this point on) \cite{CS05}. The memory terms are, in general, integral terms which account for the history of the variables that are not resolved. One would like, if possible, to compute these memory integrals {\it without} having to solve the full system. This is a difficult task, since it is rarely clear how the memory of a reduced model (which is based on the dynamics of the unresolved variables) can be estimated from pure analytical considerations or even relatively cheap numerical calculations involving only the resolved variables.  On the other hand, for problems of practical interest where the solution of the full system may be, at best, only feasible for short times, we are forced to consider ways of estimating the memory terms from such analytical or lower dimensional considerations.

We use the case of the viscous Burgers equation with uncertain viscosity coefficient to illustrate how it is possible to estimate the parameters needed to specify the memory terms. The basic idea is that the uncertainty in the viscosity coefficient leads to linear and nonlinear contributions in the memory terms. One can group the linear contributions from the different memory terms and then require that this linear term is a stabilizing one. This procedure allows to estimate recursively (as we increase the order of the terms kept in the reduced model) the parameters involved in the memory integrals.    

Section \ref{mz_formalism} presents a brief introduction to the MZ formalism for the construction of reduced models of systems of ODEs. In Section \ref{memory_comp} we develop a reformulation of the MZ formalism. This allows the  calculation of the memory terms through the solution of ordinary differential equations instead of the computation of convolution integrals as they appear in the original formulation. Section \ref{example} applies the reformulation of MZ presented in Section \ref{memory_comp} to the viscous Burgers equation when the viscosity coefficient is uncertain. Finally, in Section \ref{discussion} we discuss certain directions for future work.

\section{Mori-Zwanzig formalism}\label{mz_formalism}

We begin with a brief presentation of the Mori-Zwanzig formalism \cite{CHK00,CHK3}. Suppose we are given the system 
\begin{equation}\label{odes}
\frac{du(t)}{dt} = R (t,u(t)),
\end{equation}
where $u = ( \{u_k\}), \; k \in H \cup G$ 
with initial condition $u(0)=u_0.$ Our goal is to construct a reduced model for the modes in the subset $H.$ The system of ordinary differential equations
we are given can be transformed into a system of  linear
partial differential equations
\begin{equation}
\label{pde}
\pd{\phi_k}{t}=L \phi_k, \qquad \phi_k (u_0,0)=u_{0k}, \, k \in H \cup G
\end{equation}
where $L=\sum_{k \in H \cup G } R_i(u_0) \frac{\partial}{\partial u_{0i}}.$ The solution of \eqref{pde} is
given by $u_k (u_0,t)=\phi_k(u_0,t)$. Using semigroup notation we can rewrite (\ref{pde}) as
$$\pd{}{t} e^{tL} u_{0k}=L e^{tL} u_{0k}$$
Suppose that the vector of initial conditions can be divided as $u_0=(\hat{u}_0,\tilde{u}_0),$ where 
$\hat{u}_0$ is the vector of the resolved variables (those in $H$) and $\tilde{u}_0$ is the vector of the unresolved variables (those in $G$).  Let $P$ be an orthogonal projection on the space of functions of $\hat{u}_0$ and $Q=I-P.$ 

Equation \eqref{pde} 
can be rewritten as 
\begin{equation}
\label{mz}
\frac{\partial}{\partial{t}} e^{tL}u_{0k}=
e^{tL}PLu_{0k}+e^{tQL}QLu_{0k}+
\int_0^t e^{(t-s)L}PLe^{sQL}QLu_{0k}ds, \, k \in H,
\end{equation}
where we have used Dyson's formula
\begin{equation}
\label{dyson1}
e^{tL}=e^{tQL}+\int_0^t e^{(t-s)L}PLe^{sQL}ds.
\end{equation}
Equation (\ref{mz}) is the Mori-Zwanzig identity. 
Note that
this relation is exact and is an alternative way
of writing the original PDE. It is the starting
point of our approximations. Of course, we
have one such equation for each of the resolved
variables $u_k, k \in H$. The first term in (\ref{mz}) is
usually called Markovian since it depends only on the values of the variables
at the current instant, the second is called "noise" and the third "memory". 

If we write
$$e^{tQL}QLu_{0k}=w_k,$$ 
$w_k(u_0,t)$ satisfies the equation
\begin{equation}
\label{ortho}
\begin{cases}
&\frac{\partial}{\partial{t}}w_k(u_0,t)=QLw_k(u_0,t) \\ 
& w_k(u_0,0) = QLx_k=R_k(u_0)-(PR_k)(\hat{u_0}). 
\end{cases} 
\end{equation}
If we project (\ref{ortho}) we get
$$P\frac{\partial}{\partial{t}}w_k(u_0,t)=
PQLw_k(u_0,t)=0,$$
since $PQ=0$. Also for the initial condition
$$Pw_k(u_0,0)=PQLu_{0k}=0$$
by the same argument. Thus, the solution
of (\ref{ortho}) is at all times orthogonal
to the range of $P.$ We call
(\ref{ortho}) the orthogonal dynamics equation. Since the solutions of 
the orthogonal dynamics equation remain orthogonal to the range of $P$, 
we can project the Mori-Zwanzig equation (\ref{mz}) and find
\begin{equation}
\label{mzp}
\frac{\partial}{\partial{t}} Pe^{tL}u_{0k}=
Pe^{tL}PLu_{0k}+
P\int_0^t e^{(t-s)L}PLe^{sQL}QLu_{0k} ds.
\end{equation}

\section{Finite memory}\label{memory_comp}

In this section we describe a reformulation of the problem of computing the memory term which 
does not use the orthogonal dynamics equation. We focus on the case when the memory has a finite extent only. The case of infinite memory is simpler and is a special case of the formulation presented below. Also, the current reformulation allows us to comment on what happens in the case when the memory is very short. 

Let $w_{0k}(t)=P\int_0^t e^{(t-s)L}PLe^{sQL}QLu_{0k} ds=P\int_0^t e^{sL}PLe^{(t-s)QL}QLu_{0k} ds,$ by the change of variables $t'=t-s.$ Note, that $w_{0k}$ depends both on $t$ and the resolved part of the initial conditions $\hat{u}_0.$ We have suppressed the $\hat{u}_0$ dependence for simplicity of notation. If the memory extends only for $t_0$ units in the past (with $t_0 \leq t,$) then $$w_{0k}(t)=P\int_{t-t_0}^t e^{sL}PLe^{(t-s)QL}QLu_{0k} ds.$$ The evolution of 
$w_{0k}$ is given by 
\begin{equation}\label{memory_1}
\frac{dw_{0k}}{dt}=Pe^{tL}PLQLu_{0k}-Pe^{(t-t_0)L}PLe^{t_0 QL}QLu_{0k}+w_{1k}(t),
\end{equation}
where $$w_{1k}(t)=P\int_{t-t_0}^t e^{sL}PLe^{(t-s)QL}QLQLu_{0k} ds.$$ To allow for more flexibility, let us assume that the integrand in the formula for $w_{1k}(t)$ contributes only for $t_1$ units with $t_1 \leq t_0.$ Then $$w_{1k}(t)=P\int_{t-t_1}^t e^{sL}PLe^{(t-s)QL}QLQLu_{0k} ds.$$ 
We can proceed and write an equation for the evolution of $w_{1k}(t)$ which reads
\begin{equation}\label{memory_2}
\frac{dw_{1k}}{dt}=Pe^{tL}PLQLQLu_{0k}-Pe^{(t-t_1)L}PLe^{t_1 QL}QLQLu_{0k}+w_{2k}(t),
\end{equation}
where $$w_{2k}(t)=P\int_{t-t_1}^t e^{sL}PLe^{(t-s)QL}QLQLQLu_{0k} ds.$$ Similarly, if this integral extends only for $t_2$ units in the past with $t_2 \leq t_1,$ then
$$w_{2k}(t)=P\int_{t-t_2}^t e^{sL}PLe^{(t-s)QL}QLQLQLu_{0k} ds.$$
This hierarchy of equations continues indefinitely. Also, we can assume for more flexibility that at every level of the hierarchy we allow the interval of integration for the integral term to extend to fewer or the same units of time than the integral in the previous level. If we keep, say, $n$ terms in this hierarchy, the equation for $w_{(n-1)k}(t)$ will read 
\begin{gather}\label{memory_n}
\frac{dw_{(n-1)k}}{dt}=Pe^{tL}PL(QL)^{n-1}QLu_{0k}- \\
 Pe^{(t-t_{n-1})L}PLe^{t_{n-1} QL}(QL)^{n-1}QLu_{0k}+w_{nk}(t)  \notag
\end{gather}
where $$w_{nk}(t)=P\int_{t-t_n}^t e^{sL}PLe^{(t-s)QL}(QL)^{n}QLu_{0k} ds$$
Note that the last term in \eqref{memory_n} involves the unknown evolution operator for the orthogonal dynamics equation. This situation is the well-known closure problem. We can stop the hierarchy at the $n$th term by assuming that $w_{nk}(t)=0.$

In addition to the closure problem, the unknown evolution operator for the orthogonal dynamics equation appears in the equations for the evolution of $w_{0k}(t),\ldots,w_{(n-1)k}(t)$ through the terms $Pe^{(t-t_0)L}PLe^{t_0 QL}QLu_{0k},\ldots$ $Pe^{(t-t_0)L}PLe^{t_0 QL}(QL)^{n-1}QLu_{0k}$ respectively.

We describe now a way to express these terms involving the unknown orthogonal dynamics operator through known quantities so that we obtain a closed system for the evolution of $w_{0k}(t),\ldots,w_{(n-1)k}(t).$

Since we want to treat the case where $t_0$ is not necessarily small, we divide the interval $[t-t_0,t]$ in $n_0$ subintervals. Define 

\begin{align*}
w_{0k}^{(1)}(t) & =P\int_{t-\Delta t_0}^t e^{sL}PLe^{(t-s)QL}QLu_{0k} ds \\
w_{0k}^{(2)}(t) & =P\int_{t-2 \Delta t_0}^{t- \Delta t_0} e^{sL}PLe^{(t-s)QL}QLu_{0k} ds \\
\ldots & \\
w_{0k}^{(n_0)}(t) & =P\int_{t-t_0}^{t- (n_0-1)\Delta t_0} e^{sL}PLe^{(t-s)QL}QLu_{0k} ds,
\end{align*}
where $n_0 \Delta t_0 = t_0$ and $w_{0k}(t)=\sum_{i=1}^{n_0} w_{0k}^{(i)}(t).$ Similarly, we can define the quantities $w_{1k}^{(1)}(t),\ldots,w_{1k}^{(n_1)}(t)$ 
\begin{align*}
w_{1k}^{(1)}(t) & =P\int_{t-\Delta t_1}^t e^{sL}PLe^{(t-s)QL}QLQLu_{0k} ds \\
w_{1k}^{(2)}(t) & =P\int_{t-2 \Delta t_1}^{t- \Delta t_1} e^{sL}PLe^{(t-s)QL}QLQLu_{0k} ds \\
\ldots & \\
w_{1k}^{(n_1)}(t) & =P\int_{t-t_1}^{t- (n_1-1)\Delta t_1} e^{sL}PLe^{(t-s)QL}QLQLu_{0k} ds,
\end{align*}
where $n_1 \Delta t_1 = t_1$ and $w_{1k}(t)=\sum_{i=1}^{n_1} w_{1k}^{(i)}(t).$ In a similar fashion we can define corresponding quantities for all the memory terms up to  $w_{(n-1)k}(t)=\sum_{i=1}^{n_{n-1}} w_{(n-1)k}^{(i)}(t).$ 

In order to proceed we need to make an approximation for the integrals over the subintervals.

\subsection{Trapezoidal rule approximation}\label{trapezoidal}
We have
\begin{multline*}
w_{0k}^{(1)}(t)  =P\int_{t-\Delta t_0}^t e^{sL}PLe^{(t-s)QL}QLu_{0k} ds   \\
=\biggl[ Pe^{tL}PLQLu_{0k}+Pe^{(t-\Delta t_0)L}PLe^{\Delta t_0 QL}QLu_{0k} \biggr] \frac{\Delta t_0}{2}+ O((\Delta t_0)^3)
\end{multline*}
from which we find
$$Pe^{(t-\Delta t_0)L}PLe^{\Delta t_0 QL}QLu_{0k}=\biggl ( \frac{2}{\Delta t_0} \biggr ) w_{0k}^{(1)}(t) - Pe^{tL}PLQLu_{0k} + O((\Delta t_0)^2)$$
and from \eqref{memory_1}
\begin{equation*}
\frac{dw_{0k}^{(1)}}{dt}=-\biggl ( \frac{2}{\Delta t_0} \biggr ) w_{0k}^{(1)}(t)+ 2Pe^{tL}PLQLu_{0k}+w_{1k}^{(1)}(t)+ O((\Delta t_0)^2).
\end{equation*}
Similarly, for $w_{0k}^{(2)}(t)$ we find
\begin{multline*}
\frac{dw_{0k}^{(2)}}{dt}=\biggl ( \frac{4}{\Delta t_0} \biggr ) w_{0k}^{(1)}(t) \\
-\biggl ( \frac{2}{\Delta t_0} \biggr ) w_{0k}^{(2)}(t) - 2Pe^{tL}PLQLu_{0k} 
+w_{1k}^{(2)}(t)+ O((\Delta t_0)^2)
\end{multline*}
In general,
\begin{multline}\label{memory_1a}
\frac{dw_{0k}^{(i)}}{dt}= -\biggl ( \frac{2}{\Delta t_0} \biggr ) w_{0k}^{(i)}(t) + (-1)^{i+1} 2Pe^{tL}PLQLu_{0k} \\
 +\biggl [  \sum_{j=1}^{i-1}  \biggl ( \frac{4}{\Delta t_0} \biggr ) (-1)^{i+j+1} w_{0k}^{(j)}(t) \biggr ] +w_{1k}^{(i)}(t)+ O((\Delta t_0)^2)  \; \; \text{for}  \; \; i=1,\ldots,n_0.
\end{multline}
Similarly,
\begin{multline*}
\frac{dw_{1k}^{(i)}}{dt}= -\biggl ( \frac{2}{\Delta t_1} \biggr ) w_{1k}^{(i)}(t) + (-1)^{i+1} 2Pe^{tL}PLQLQLu_{0k} \\
 +\biggl [  \sum_{j=1}^{i-1}  \biggl ( \frac{4}{\Delta t_1} \biggr ) (-1)^{i+j+1} w_{1k}^{(j)}(t) \biggr ] +w_{2k}^{(i)}(t)+ O((\Delta t_1)^2)  \; \; \text{for}  \; \; i=1,\ldots,n_1 
\end{multline*} 
$\ldots$
\begin{multline}
\frac{dw_{(n-1)k}^{(i)}}{dt}= -\biggl ( \frac{2}{\Delta t_{n-1}} \biggr ) w_{(n-1)k}^{(i)}(t) + (-1)^{i+1} 2Pe^{tL}PL(QL)^{n-1}QLu_{0k} \\
 +\biggl [  \sum_{j=1}^{i-1}  \biggl ( \frac{4}{\Delta t_{n-1}} \biggr ) (-1)^{i+j+1} w_{(n-1)k}^{(j)}(t) \biggr ] + O((\Delta t_{n-1})^2)  \; \; \text{for}  \; \; i=1,\ldots,n_{n-1}.
\end{multline}
By dropping the $O((\Delta t_0)^2),\ldots, O((\Delta t_{n-1})^2)$ terms we obtain a system of $n_0+n_1+\ldots+n_{n-1}$ differential equations for the evolution of the quantities $w_{0k}^{(1)}(t),\ldots,w_{(n-1)k}^{(n_{n-1})}.$ This system allows us to determine the memory term $w_{0k}(t)=P\int_0^t e^{(t-s)L}PLe^{sQL}QLu_{0k} ds.$ Since the approximation we have used for the integral leads to an error $O(\Delta t)^2,$ the ODE solver should also be $O(\Delta t)^2.$ We have used the modified Euler method to solve numerically the equations for the reduced model. 

Note that the implementation of the above scheme requires the knowledge of the expressions for $Pe^{tL}PLQLu_{0k},\ldots,Pe^{tL}PL(QL)^{n-1}QLu_{0k}.$ Since the computation of these expressions for large $n$ can be rather involved for nonlinear systems (see Section \ref{example}), we expect that the above scheme will be used with a small to moderate value of $n.$ Finally, we mention that the above construction can be carried out for integration rules of higher order e.g. Simpson's rule. 

\section{Burgers equation with uncertain viscosity coefficient}\label{example}

In this section we show how the above MZ formulation can be used for uncertainty quantification (UQ). In particular, we apply it to the one-dimensional Burgers equation when the viscosity coefficient is uncertain. The equation is given by 
\begin{equation}\label{burgersequation}
u_t+u u_x = \nu u_{xx},
\end{equation}
where $\nu > 0.$ Equation (\ref{burgersequation}) should be supplemented with an initial condition $u(x,0)=u_0(x)$ and boundary conditions. We solve (\ref{burgersequation}) in the interval $[0,2\pi]$ with periodic boundary conditions. This allows us to expand the solution in Fourier series
$$u_{N}(x,t )=\underset{k \in F}{\sum} u_k(t) e^{ikx},$$
where $F=[-\frac{N}{2},\frac{N}{2}-1].$ The equation of motion for the Fourier mode $u_k$ becomes
\begin{equation}
\label{burgersode}
 \frac{d u_k}{dt}=- \frac{ik}{2} \underset{p, q \in F}{\underset{p+q=k  }{ \sum}} u_{p} u_{q}  -\nu k^2 u_k.
\end{equation}
We assume that the viscosity coefficient $\nu$ is uncertain (random) and can be expanded as $\nu (\xi)= \nu_0+ \alpha \xi$ where $\xi$ is uniformly distributed in $[-1,1].$ In the numerical experiments we have taken $\nu_0=0.1$ and $\alpha=0.07.$ This means that the viscosity coefficient is allowed to take values in the interval $[.03,1.07].$ The choice of the range allows us to compute an accurate solution for any viscosity coefficient in the range without having to employ a large number of Fourier modes. 

To proceed we expand the solution $u_k(t,\xi)$ for $k \in F$ in a polynomial chaos expansion using Legendre polynomials which are orthogonal in the interval $[-1,1].$ In particular, we have that $$\int_{-1}^1L_i (\xi) L_j(\xi) d \xi=\frac{2}{2i+1} \delta_{ij},$$ where $L_i(\xi)$ is the Legendre polynomial of order $i.$ For each wavenumber $k$ we expand the solution $u_k(t,\xi)$ of \eqref{burgersodemz} in Legendre polynomials and keep the first $M$ polynomials  
\begin{equation}\label{ode_expansion}
u_k(t,\xi)\approx \sum_{i=0}^{M-1} u_{ki}(t) L_i(\xi), \; \; \text{where} \; \; \xi \sim U[-1,1].
\end{equation}
Similarly, the viscosity coefficient can be written as $\nu = \sum_{i=0}^1\nu_i L_i(\xi)$ with $\nu_1=\alpha$  since $L_0 (\xi)=1$ and $L_1 (\xi)=\xi.$
Substitution of \eqref{ode_expansion} in \eqref{burgersodemz}, use of the expansion of the viscosity coefficient and of the orthogonality properties of the Legendre polynomials gives 
\begin{equation}\label{burgersodemz_system}
\frac{du_{kr}(t)}{dt}=- \frac{ik}{2} \sum_{l=0}^{M-1} \sum_{m=0}^{M-1} \underset{p, q \in F}{\underset{p+q=k  }{ \sum}} u_{pl} u_{qm}  c_{lmr} - k^2 \sum_{l=0}^{M-1} \sum_{m=0}^{M-1}\nu_l  u_{km} c_{lmr}
\end{equation}
for $k \in F$ and $r=0,\ldots,M-1.$ Also $$c_{lmr}=\frac{E[ L_l (\xi) L_m(\xi) L_r(\xi) ] }{E[L^2_r(\xi) ]},$$
where the expectation $E[\cdot]$ is taken with respect to the uniform density on $[-1,1].$ The expectation on the denominator of the expression for $c_{lmr}$ is $E[L^2_r(\xi) ]=\int_{-1}^1 L^2_r(\xi) \frac{1}{2}d\xi=\frac{1}{2r+1},$ while the expectation on the numerator can be computed accurately using Gaussian quadrature with Legendre nodes. The Legendre polynomial triple product integral defines a tensor which has the following sparsity pattern: $E[ L_l (\xi) L_m(\xi) L_r(\xi) ]=0,$ if $ l+m < r$ or $l+r < m$ or $m+r < l$ or $l+m+r= \text{odd}$ \cite{gupta}. Due to this sparsity pattern, for a given value of $M$ only about $1/4$ of the $M^3$ tensor entries are different from zero. The sparsity pattern will be used below (see Section \ref{mz_memory}) to facilitate the estimation of the length of the memory.

\subsection{MZ reduced model}\label{mz_ode_example}

To conform with 
the Mori-Zwanzig formalism we set 
$$R_{kr}(u)=- \frac{ik}{2} \sum_{l=0}^{M-1} \sum_{m=0}^{M-1} \underset{p, q \in F}{\underset{p+q=k  }{ \sum}} u_{pl} u_{qm}  c_{lmr} -k^2 \sum_{l=0}^{M-1} \sum_{m=0}^{M-1}\nu_l  u_{km} c_{lmr},$$
where $u=\{u_{kr}\}$ for $k \in F$ and $r=0,\ldots,M-1.$ Thus, we have
\begin{equation}
\label{burgersodemz}
\frac{d u_{kr}}{dt}=R_{kr}(u) 
\end{equation}
for $k \in F$ and $r=0,\ldots,M-1.$  
We proceed by dividing the variables in resolved and unresolved. In particular, we consider as resolved the variables $\hat{u}=\{u_{kr}\}$ for $k \in F$ and $r=0,\ldots,\Lambda-1,$ where $\Lambda < M.$ Similarly, the unresolved variables are $\tilde{u}=\{u_{kr}\}$ for $k \in F$ and $r=\Lambda,\ldots,M-1.$ In the notation of Section \ref{mz_formalism} we have $H= F \cup (0,\ldots,\Lambda-1)$ and $G= F\cup (\Lambda,\ldots,M-1).$ In other words, we resolve, for all the Fourier modes, only the first $\Lambda$ of the Legendre expansion coefficients and we shall construct a reduced model for them. 

The system (\ref{burgersodemz}) is supplemented by the initial 
condition $u_0=(\hat{u}_0,\tilde{u}_0).$ We focus on initial conditions where 
the unresolved Fourier modes are set to zero, i.e. $u_0=(\hat{u}_0,0).$ We also define $L$ by 
$$L=\sum_{k \in F}\sum_{r=0}^{M-1} R_{kr}(u_0) \frac{\partial}{\partial u_{0kr}}.$$ 
To construct a MZ reduced model we need to define a projection operator $P.$ For a function $h(u_0)$ of all the 
variables, the projection operator we will use is defined by $P(h(u))=P(h(\hat{u}_0,\tilde{u}_0))=h(\hat{u}_0,0),$ i.e. 
it replaces the value of the unresolved variables $\tilde{u}_0$ in any function $h(u_0)$ by zero. Note that this choice of projection is consistent with the initial conditions we have chosen. Also, we define the Markovian term 
$$ PLu_{0k}=PR_k(u_0)=- \frac{ik}{2} \sum_{l=0}^{\Lambda-1} \sum_{m=0}^{\Lambda-1} \underset{p, q \in F}{\underset{p+q=k  }{ \sum}} u_{0pl} u_{0qm}  c_{lmr} -k^2 \sum_{l=0}^{M-1} \sum_{m=0}^{\Lambda-1}\nu_l  u_{0km} c_{lmr}.$$ 
The Markovian term has the same functional form as the RHS of the full system but is restricted to a sum over only the first $\Lambda$ Legendre expansion coefficients  for each Fourier mode. 

\subsection{Number of memory terms and memory length}\label{mz_memory}

\subsubsection{Number of memory terms}\label{mz_number}
We have to decide on the number of terms that will be used in the expansion of the memory as well as the length of the memory kept for each term (see Section \ref{memory_comp}). The fact that we are considering the case of uncertain viscosity coefficient becomes important in choosing how many terms to keep in the memory expansion and what the memory length should be for each term. To see this we need to compute the first few terms in the expansion. For the first two terms $PLQLu_{0kr}$ and $PLQLQLu_{0kr}$ we find

\begin{gather}\label{burgersmemory1}
PLQLu_{0kr}=2\times \biggl [   - \frac{ik}{2}   \sum_{l=\Lambda}^{M-1} \sum_{m=0}^{\Lambda-1} \underset{p, q \in F}{\underset{p+q=k  }{ \sum}} PLu_{0pl} u_{0qm}  c_{lmr} \biggr ]  \\
-k^2 \sum_{l=0}^{M-1} \sum_{m=\Lambda}^{M-1}\nu_l  PLu_{0km} c_{lmr} \notag
\end{gather}
and
\begin{gather}\label{burgersmemory2}
PLQLQLu_{0kr}=2\times \biggl [   - \frac{ik}{2}   \sum_{l=\Lambda}^{M-1} \sum_{m=0}^{\Lambda-1} \underset{p, q \in F}{\underset{p+q=k  }{ \sum}} PLQLu_{0pl} u_{0qm}  c_{lmr} \biggr ]  \\
+2\times \biggl [   - \frac{ik}{2}   \sum_{l=\Lambda}^{M-1} \sum_{m=0}^{M-1} \underset{p, q \in F}{\underset{p+q=k  }{ \sum}} PLu_{0pl} PLu_{0qm}  c_{lmr} \biggr ]  \notag \\
-k^2 \sum_{l=0}^{M-1} \sum_{m=\Lambda}^{M-1}\nu_l  PLQLu_{0km} c_{lmr} \notag
\end{gather}
For the sake of simplicity, we restrict attention to the case when $\Lambda=1,$ so that we resolve only the zeroth term in the Legendre expansion. The linear (viscous) term in \eqref{burgersodemz_system} contributes a linear destabilizing term in $PLQLu_{0k0}$ and a linear stabilizing term in $PLQLQLu_{0k0}.$ To show this, we use the fact that the Legendre expansion of the viscosity coefficient has only the zero and first components nonzero, the properties of the defined projection operator $P$ and the sparsity of the Legendre polynomial triple product. Through straightforward but tedious algebra we find that the viscous term in \eqref{burgersodemz_system} contributes the term $$k^4 \nu_1^2  c_{101} c_{110} u_{0k0}$$ in $PLQLu_{0k0}$ and the term $$-k^6 \nu_0 \nu_1^2 c_{011}c_{101} c_{110} u_{0k0}$$ in $PLQLQLu_{0k0}.$ Indeed, the contribution to $PLQLu_{0k0}$ is destabilizing and the contribution to $PLQLQLu_{0k0}$ is a stabilizing term. Note that these contributions correspond to terms of the form $u_{xxxx}$ and $u_{xxxxxx}$ respectively in real space. 

With more effort one can compute the contributions of the linear viscous term to the memory terms $PLQLQLQLu_{0k0}$ and $PLQLQLQLQLu_{0k0}.$ One finds that the viscous term contributes the term $$k^8 [\nu_0^2 \nu_1^2 c_{110} c_{011}^2 c_{101}+ \nu_1^4 c_{110} c_{121} c_{101} c_{112}]  u_{0k0}$$ to $PLQLQLQLu_{0k0}$ and the term  
\begin{gather*}
- k^{10} [\nu_0^3 \nu_1^2 c_{110} c_{011}^3 c_{101}+ \nu_0 \nu_1^4 c_{110} c_{011} c_{121} c_{101} c_{112} \\
 + \nu_0 \nu_1^4 c_{110} c_{121} c_{022} c_{101} c_{112}+ \nu_0 \nu_1^4 c_{110}c_{121}c_{112}c_{011}c_{101}] u_{0k0}
 \end{gather*}
  to $PLQLQLQLQLu_{0k0}.$ 
  
The pattern of alternating destabilizing and stabilizing contributions of the viscous term to the memory terms continues for higher order terms. The conclusion from this pattern is that one needs to keep the memory terms in pairs in order to guarantee the stability of the reduced model. Also, since these contributions correspond to higher and higher spatial derivatives in real space we have to use only a few pairs otherwise the reduced model will become extremely stiff. In our numerical experiments we have used only the first pair of memory terms, namely $PLQLu_{0k0}$ and $PLQLQLu_{0k0}.$

\subsubsection{Length of the memory}\label{mz_length}
We will use the information obtained in Section \ref{mz_number} to estimate the length of the memory. Ideally, we would like to estimate the values of $t_0$ and $t_1$ {\it without} having to solve the full system. That would make the construction of the reduced model efficient and applicable in cases where the solution of the full system is expensive (or possibly unknown). 

We focus again on the case when $\Lambda=1$ and we assume that we use only one subinterval to discretize the time integrals, i.e. $\Delta t_0=t_0$ and $\Delta t_1=t_1.$ If we keep only the terms $PLQLu_{0k0}$ and $PLQLQLu_{0k0}$ for the memory, the reduced model reads

\begin{gather}
\frac{d u_{k0}}{dt}= e^{tL}PLu_{0k0} + w_{0k0}(t) \label{reduced1} \\
\frac{d w_{0k0}}{dt}=2Pe^{tL}PLQLu_{0k0} - \frac{2}{t_0} w_{0k0}(t) + w_{1k0}(t) \label{reduced2}\\
\frac{d w_{1k0}}{dt}=2Pe^{tL}PLQLQLu_{0k0} - \frac{2}{t_1} w_{1k0}(t) \label{reduced3}
\end{gather}
We can solve \eqref{reduced2} and \eqref{reduced3} formally and substitute in \eqref{reduced1} to get
\begin{gather}\label{reduced_integral}
\frac{d u_{k0}}{dt}= e^{tL}PLu_{0k0} + \int_0^t e^{-\lambda_0(t-s)}2Pe^{sL}PLQLu_{0k0} ds \\
+ \int_0^t e^{-\lambda_0(t-s)}  \int_0^s e^{-\lambda_1(s-\tau)}2Pe^{\tau L}PLQLQLu_{0k0} d\tau ds,  \notag
\end{gather}
where $\lambda_0=2/t_0$ and $\lambda_1=2/t_1.$ 

The quantities  $w_{0k0}(t)$ and $w_{1k0}(t)$ have on the RHS of their equations of evolution the terms $Pe^{tL}PLQLu_{0k0}$ and $Pe^{tL}PLQLQLu_{0k0}$ respectively. As we have seen, the linear contributions of the viscous term to $PLQLu_{0k0}$ and $PLQLQLu_{0k0}$  correspond to higher spatial derivatives. Due to the presence of $k^4$ and $-k^6$ in the linear terms in the expressions for $PLQLu_{0k0}$ and $PLQLQLu_{0k0}$ respectively, those linear terms are going to have large value for large wavenumbers. At the same time, these linear terms are linear in $u_{k0}(t)$ which is expected to evolve more slowly. Thus, the linear terms in the expressions for $PLQLu_{0k0}$ and $PLQLQLu_{0k0}$ are expected to have (at least for large wavenumbers) large values and evolve slowly. As a result, we expect the quantities  $w_{0k0}(t)$ and $w_{1k0}(t)$ to evolve faster than $u_{k0}(t).$ With this in mind, we expect the memory lengths $t_0$ and $t_1$ to be shorter compared to the time scale of evolution of $u_{k0}(t).$ The crudest approximation that one can make for the integrals in \eqref{reduced_integral} are $$\int_0^t e^{-\lambda_0(t-s)}2e^{sL}PLQLu_{0k0} ds \approx \frac{1}{\lambda_0} 2Pe^{tL}PLQLu_{0k0} = \frac{t_0}{2} 2Pe^{tL}PLQLu_{0k0} $$ and $$ \int_0^t e^{-\lambda_0(t-s)}  \int_0^s e^{-\lambda_1(s-\tau)}2e^{\tau L}PLQLQLu_{0k0}  \approx  \frac{t_0 t_1}{4} 2Pe^{tL}PLQLQLu_{0k0}.$$
These approximations for the integrals allow us to group together all the linear terms on the RHS of the equation for $u_{k0}(t).$ Indeed, putting together the linear contributions from the Markovian term and the two integral terms (after the approximation) we get the linear term
\begin{equation}\label{stability}
\biggl [ -k^2 \nu_0 c_{000} + t_0 k^4 \nu_1^2  c_{101} c_{110} - \frac{t_0 t_1}{2}k^6 \nu_0 \nu_1^2 c_{011}c_{101} c_{110} \biggr ] u_{k0}(t).
\end{equation}
The expression in brackets in \eqref{stability} can be used to determine what should the values of $t_0$ and $t_1$ be so that the reduced model is linearly stable for all wave numbers $k \in F.$ 

The calculation of $t_0$ and $t_1$ is done in two steps. If we ignore the $k^4$ and $k^6$ terms, then the only contribution is from the Markovian term and the bracketed expression is a parabola in $k$ with negative values. If we also include the $k^4$ term the parabola changes into a double-well curve which can become greater than zero for some wave numbers depending on the value of $t_0.$ In fact, we can estimate the minimum value of $t_0$ for this to happen by solving the equation  $$-k_{max}^2 \nu_0 c_{000} + t^{min}_0 k_{max}^4 \nu_1^2 c_{101} c_{110} =0$$ where $k_{max}$ is the maximum wavenumber present in the solution. In fact, $t^{min}_0= (\nu_0 c_{000} )/( k_{max}^2\nu_1^2 c_{101} c_{110}).$ For $t^{min}_0$ all the wavenumbers are linearly stable except for the wavenumber $k_{max}$ which is only marginally stable. 

Now, suppose that we set $t_0$ equal to $t^{min}_0$ and we also include the $k^6$ term. Then, the bracketed expression in \eqref{stability} becomes a 6th order negative curve. The addition of the negative definite $k^6$ term provides us with an advantage. It allows us to increase $t_0$ to values larger than $t^{min}_0$ as long as $t_1$ is appropriately chosen to make sure that the all the wavenumbers are linearly stable. Of course, one should not increase $t_0$ too much because a correspondingly large value of $t_1$ in conjunction with the $k^6$ factor can render the reduced model very stiff. Thus, the final criterion which allows us to determine $t_0$ and $t_1$ uniquely is that, based on the linear stability domain of the numerical method, we pick $t_0$ and $t_1$ so that the required step size for the reduced model is not smaller than the step size for the original (full) system.      

\subsection{Numerical results}\label{numerical}

In this section we present numerical results for the reduced model of the viscous Burgers equation with uncertain viscosity coefficient given by $\nu = \sum_{i=0}^1\nu_i L_i(\xi)$ with $\nu_0=0.1$ and $\nu_1=0.07.$ The solution of the full system was computed with $N=96$ Fourier modes ($F=[-48,47]$) and the first 7 Legendre polynomials ($M=7$). The first 7 Legendre polynomials were enough to obtain converged statistics for the full system. The full system was solved with the modified Euler method with $\Delta t = 0.001.$

The reduced model uses $N=96$ Fourier modes but only the first Legendre polynomial, so $\Lambda=1.$ The memory length parameters in the reduced model were chosen to be $t_0=0.2$ and $t_1=0.01632$ according to the scheme presented in Section \ref{mz_length}. In particular, this choice of memory length guarantees linear stability of the reduced model when it is solved with the modified Euler method with a step size of $\Delta t =0.001.$ We discretize the memory integral with 1 subinterval, i.e. $n_0=n_1=1,$ $\Delta t_0=0.2$ and $\Delta t_1=0.01632$ according to the notation of Section \ref{memory_comp}. With this choice of parameters the running time of the reduced model is about half of that of the full system. 

Note that according to our analysis in Section \ref{memory_comp} there is a discrepancy in the local truncation error estimates of the trapezoidal rule and the modified Euler scheme. For the trapezoidal rule the local truncation error estimate is $O((0.2)^2)$ for the memory term $w_{0k}$ and $O((0.01632)^2)$ for the memory term $w_{1k}.$ On the other hand, the local truncation error estimate for the modified Euler method is $O((0.001)^2).$ To make the discrepancy disappear we must use more subintervals at the cost of making the reduced model evolution more expensive. We tried that but the accuracy of the results of the reduced model did not change.

\begin{figure}
\centering
\epsfig{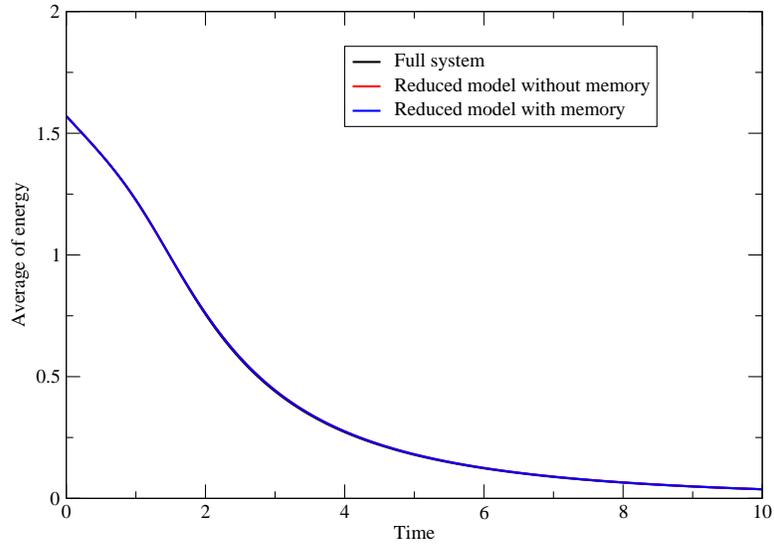}
\caption{Evolution of the mean of the energy of the solution using only the first Legendre polynomial.}
\label{plot_viscosity_energy}
\end{figure}

Figure \ref{plot_viscosity_energy} shows the evolution of the mean energy of the solution 
$$E(t)=\frac{1}{2}\sum_{k \in F} 2\pi |u_{k0}(t)|^2$$
as computed from the full system (with $M=7$ Legendre polynomials), the MZ reduced model with $\Lambda=1$ {\it without} memory (keeping only the Markovian term) and the MZ reduced model with $\Lambda=1$ {\it with} memory. The reduced model performs equally well with or without memory.

\begin{figure}
\centering
\epsfig{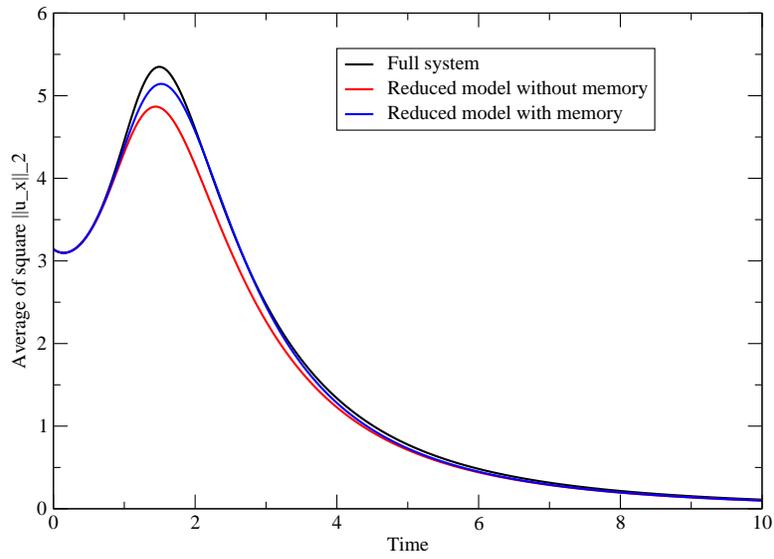}
\caption{Evolution of the mean of the squared $l_2$ norm of the gradient of the solution calculated using only the first Legendre polynomial.}
\label{plot_viscosity_gradient}
\end{figure}

Figure \ref{plot_viscosity_gradient} shows the evolution of the mean squared $l_2$ norm of the gradient of the solution 
$$G(t)=\sum_{k \in F} 2 \pi k^2 |u_{k0}(t)|^2$$
as computed from the full system (with $M=7$ Legendre polynomials), the MZ reduced model with $\Lambda=1$ {\it without} memory (keeping only the Markovian term) and the MZ reduced model with $\Lambda=1$ {\it with} memory. It is obvious from the figures that the inclusion of the memory term improves considerably the performance of the reduced model. 

By looking at Figure \ref{plot_viscosity_gradient}, we see that the reduced model with memory predicts a smaller value for the peak of $G(t).$ We know that the term $PLQLu_{0k0}$ contributes a linear destabilizing term to the reduced model (see \eqref{stability}). So,  the obvious question to ask is if one can improve the accuracy of the reduced model with memory by increasing $t_0$ and correspondingly $t_1$ at the cost of making the evolution of the reduced model more expensive. As explained in Section \ref{mz_length}, the increase in the cost will come from the increased stiffness of the reduced model. We have tried increasing $t_0$ and $t_1$ and the results did not become more accurate. The reason for this lack of improvement is a sign that if one wishes to improve the accuracy of the reduced model, one needs to include higher order terms. In particular, one will have to add at least the next pair of terms, namely $PLQLQLQLu_{0k0}$ and  $PLQLQLQLQLu_{0k0}$ (see discussion at the end of Section \ref{mz_number} as to why the terms need to be added in pairs).

\section{Discussion and future work}\label{discussion}

We have presented the application of the Mori-Zwanzig formalism to the construction of reduced models for systems of differential equations resulting from polynomial chaos expansions of solutions of differential equations with parametric uncertainty. In particular, we presented a way that the reduced model can be reformulated so that instead of integro-differential equations one has to solve differential equations. The problem that arises in any reduced model with memory is to compute the length of the memory. If possible, one wishes to obtain the length of the memory without having to solve the full system. For the case of the viscous Burgers equation with uncertain viscosity coefficient, we presented a way to actually compute the length of the memory {without} having to solve the full system. Note that this construction readily applies also to other equations that include viscous dissipation e.g. the Navier-Stokes equations. 

The viscous Burgers example highlights two important issues that arise when one wants to construct reduced models for parametric uncertainty quantification. 

The first issue is how to estimate the length of the memory integrals when we keep in the reduced model more than one coefficient in the Legendre expansion, so that $\Lambda > 1.$ In this case, the bracketed expression in \eqref{stability} will be replaced by a $\Lambda \times \Lambda$ matrix. In order for the reduced model to be linearly stable, we have to require that the matrix is negative definite (or at least semidefinite). Since the elements of the matrix depend on the quantities $t_0,t_1,\ldots,$ we can use the negative definite restriction to estimate $t_0,t_1,\ldots.$

The second issue is also related to the length of the memory but addresses a different aspect. We have seen for viscous Burgers that because the memory terms correspond to higher derivatives in physical space, the lengths of the integrals for the different memory terms (in our example $t_0$ and $t_1$) {\it decrease} as the order of the memory term increases.  This allowed us to use a crude short-memory approximation of the memory integral (see Section \ref{mz_length}). On the other hand, there are cases when the length of the integrals for the different memory terms can {\it increase} as the order of the memory term increases. In such cases the short-memory approximation of the memory integrals will not work. 

A simple example which illustrates this behavior is that of a single decaying linear ode $$\frac{du}{dt}=-ku$$ where $k \sim U[0,1].$ If one applies the procedure outlined in Sections \ref{memory_comp} and \ref{trapezoidal}, it is easy to see that the terms $PLQLu_{0k0}, PLQLQLu_{0k0}, \ldots$ decrease in amplitude and thus the required memory integrals lengths $t_0,t_1, \ldots$ increase as we go to higher order terms. A crude approximation is to assume that all the memory kernels in \eqref{reduced_integral} (and for higher order terms) have the same decaying characteristic times, that is $t_0=t_1=\ldots).$ Preliminary numerical calculations show that as we increase the order of the memory terms kept in the reduced model we also have to {\it increase} the value of $t_0$ in order to increase the accuracy. A detailed analysis will be presented elsewhere.


When the uncertainty is due not to a parameter in the equations but due to the initial conditions, the criterion presented in Section \ref{mz_length} for selecting the length of the memory will not work. For example, in the viscous Burgers equation, where the viscous term is diagonal in Fourier space, if the viscosity has no uncertainty, the projection operator makes the viscous term part of the Markovian term. As a result, the viscous term will end up contributing in the memory terms but the corresponding contribution is a term which is {\it nonlinear} in the resolved variables. Thus, even if one groups the memory contributions from the viscous term, there is no simple linear stability criterion, like the one invoked in the current work, to facilitate the estimation of the memory length. The construction of reduced models for the case of uncertain initial conditions will be presented in a forthcoming publication \cite{s13}. Such a construction can also be applied to the problem of constructing reduced models for systems forced by random noise \cite{hou}. 

Finally, we mention that one can construct models which effect reduction both for the variables needed to describe uncertainty and the number of variables needed to describe the system for one realization of the uncertainty sources. This two-level reduction is imperative in situations where solving even for one realization of the uncertainty sources is very expensive e.g. atmospheric flows, fluid structure interactions.

\section*{Acknowledgements}
I would like to thank Prof. G. Karniadakis and Dr. D. Venturi for useful discussions and comments.

\end{document}